\documentclass[12pt]{article}
\usepackage[a4paper,margin=1.1in]{geometry}
\usepackage{amssymb}
\usepackage{amsmath}
\usepackage{amsfonts}
\usepackage{amsthm}
\usepackage{color}
\usepackage{float}
\usepackage[all]{xy}
\usepackage{scalerel}
\usepackage{mathabx}
\usepackage{delimset}
\usepackage{setspace}

\newcommand{\qbinom}{\genfrac{[}{]}{0pt}{}}

\newcommand{\minus}{\textup{\texttt{-}}}

\def\E{\mathop{\mbox{\textup{E}}}\nolimits}
\def\P{\mathop{\mbox{\textup{P}}}\nolimits}

\newcommand{\C}{\mathbb{C}}

\newcommand{\N}{\mathbb{N}}

\newcommand{\e}{\mathrm{e}}

\newcommand{\s}{\mathrm{S}}
\newcommand{\rr}{\mathrm{R}}

\newtheorem{theorem}{Theorem}
\newtheorem{lemma}{Lemma}
\newtheorem{definition}{Definition}
\newtheorem{proposition}{Proposition}
\newtheorem{corollary}{Corollary}

\begin{document}

\title{\textbf{Generalizing the Index of the $u$-Deformed Homogeneous Polynomials and Generating Functions for $\rr_{\minus n}(x,y;q|u)$}}
\author{Ronald Orozco L\'opez}
\newcommand{\Addresses}{{
  \bigskip
  \footnotesize

  \textit{E-mail address}, R.~Orozco: \texttt{rj.orozco@uniandes.edu.co}
  
}}

\maketitle
\tableofcontents

\begin{abstract}
This paper introduces the $u$-deformed homogeneous functions $\mathrm{R}_{\alpha}(x,y;u|q)$, for all $\alpha\in\mathbb{C}$. Basic properties of the functions $\mathrm{R}_{\alpha}(x,y;u|q)$ are given, along with recurrence relations, their $q$-difference equation, and representations. Generating functions for the functions $\mathrm{R}_{\minus n}(x,y;u|q)$ via the $u$-deformed $q$-exponential operator $\E(qD_{q}|u)$, when $u=q^b,q^{b+1/2}$, $b\geq1$, are obtained. This allows us to obtain some $q^b$-series and $q^{b+1/2}$-series. Additionally, transformation formulas for basic hypergeometric series $_{1}\phi_{1}$ and $_{1}\phi_{2}$  are derived.
\end{abstract}


\section{Introduction}

In \cite{orozco}, $u$-deformed homogeneous polynomials $\rr_{n}(x,y;u|q)$ were introduced to unify classical polynomials, such as Rogers-Szeg\"o polynomials, $q$-shifted polynomials, and Stieltjes-Wigert polynomials, among others. These polynomials are of great importance due to their applications in $q$-series identities, quantum mechanics, orthogonal polynomials theory, and combinatorics. In this paper, we generalize the polynomials $\rr_{n}(x,y;u|q)$ by extending the index $n$ to all complex numbers. That is, we define the $u$-deformed homogeneous function of order $\alpha$, for $\alpha\in\C$, as
\begin{equation}
    \rr_{\alpha}(x,y;u|q)=\sum_{k=0}^{\infty}\qbinom{\alpha}{k}_{q}u^{\binom{k}{2}}x^{\alpha\minus k}y^{k}.
\end{equation}
Using the $u$-deformed $q$-exponential operator
\begin{equation}
    \E(yD_{q}|u)=\sum_{n=0}^{\infty}u^{\binom{n}{2}}\frac{(yD_{q})^n}{(q;q)_{n}}
\end{equation}
we obtain the representation
\begin{equation}
    \rr_{\alpha}(x,y;u|q)=\E(yD_{q}|u)\{x^\alpha\}.
\end{equation}
Through the operator $\E(yD_{q}|u)$ we will find $q$-series identities for the functions 
\begin{equation}
    \rr_{-n}(x,y;u|q)=x^{-n}\sum_{k=0}^{\infty}\qbinom{n+k-1}{k}_{q}(u/q)^{\binom{k}{2}}(-y/q^nx)^k
\end{equation}
when the parameter $u$ takes the values $q^b$ and $q^{b+1/2}$, $b\geq1$. Generating functions of the functions $\rr_{\minus n}(x,y;u|q)$ are related to the $u$-deformed $q$-exponential function $\e_{q}(a/x,u)$ when the operator $\E(yD_{q}|u)$ acting on it.

\section{Preliminaries}

Some notation and terminology for basic hypergeometric series are taken from \cite{gasper}. Let $\vert q\vert<1$ and the $q$-shifted factorial be defined by
\begin{align*}
    (a;q)_{n}&=\prod_{k=0}^{n\minus1}(1-q^{k}a),\\
    (a;q)_{\infty}&=\lim_{n\rightarrow\infty}(a;q)_{n}=\prod_{k=0}^{\infty}(1-aq^{k}).
\end{align*}
The multiple $q$-shifted factorials are defined by
\begin{align*}
    (a_{1},a_{2},\ldots,a_{m};q)_{n}&=(a_{1};q)_{n}(a_{2};q)_{n}\cdots(a_{m};q)_{n},\\
    (a_{1},a_{2},\ldots,a_{m};q)_{\infty}&=(a_{1};q)_{\infty}(a_{2};q)_{\infty}\cdots(a_{m};q)_{\infty}.
\end{align*}
Some useful identities for $q$-shifted factorial:
\begin{align}
    (a;q)_{n}&=\sum_{k=0}^{n}\qbinom{n}{k}_{q}(\minus a)^{k}q^{\binom{k}{2}},\label{iden3}\\
    (a;q)_{n}&=\frac{(a;q)_{\infty}}{(aq^n;q)_{\infty}},\label{iden4}\\
    (a;q)_{n+k}&=(a;q)_{n}(aq^{n};q)_{k},\label{iden5}\\
    (aq^{\minus n};q)_{n}&=(\minus a)^nq^{\minus\binom{n+1}{2}}(a^{\minus1}q;q)_{n}.\label{iden6}
\end{align}
The $q$-binomial coefficient is defined by
\begin{equation*}
\qbinom{n}{k}_{q}=\frac{(q;q)_{n}}{(q;q)_{k}(q;q)_{n\minus k}}.
\end{equation*}
The $q$-binomial coefficient verifies that:
\begin{align}
    \qbinom{\alpha+1}{k}_{q}&=\qbinom{\alpha}{k}_{q}+q^{\alpha+1\minus k}\qbinom{\alpha}{k\minus1}_{q}=
    q^{k}\qbinom{\alpha}{k}_{q}+\qbinom{\alpha}{k\minus1}_{q},\label{eqn_pascal}\\
    \qbinom{\alpha}{k}_{q}&=\frac{(q^{\minus\alpha};q)_{k}}{(q;q)_{k}}(\minus q^\alpha)^{k}q^{\minus\binom{k}{2}},\label{iden1}\\
    \qbinom{\minus\alpha}{k}_{q}&=\qbinom{\alpha+k\minus1}{k}_{q}(\minus q^{\minus\alpha})^kq^{\minus\binom{k}{2}}.\label{iden2}
\end{align}
In our work, we will use the identities for binomial coefficients:
\begin{align}
    \binom{n+k}{2}&=\binom{n}{2}+\binom{k}{2}+nk,\label{iden7}\\
    \binom{n-k}{2}&=\binom{n}{2}+\binom{k}{2}+k(1-n).\label{iden8}
\end{align}
The ${}_r\phi_{s}$ basic hypergeometric series is define by
\begin{equation*}
    {}_r\phi_{s}\left(
    \begin{array}{c}
         a_{1},a_{2},\ldots,a_{r} \\
         b_{1},\ldots,b_{s}
    \end{array}
    ;q,z
    \right)=\sum_{n=0}^{\infty}\frac{(a_{1},a_{2},\ldots,a_{r};q)_{n}}{(q;q)_{n}(b_{1},b_{2},\ldots,b_{s};q)_{n}}\Big[(\minus1)^{n}q^{\binom{n}{2}}\Big]^{1+s-r}z^n.
\end{equation*}
In this paper, we will frequently use the $q$-binomial theorem:
\begin{equation}
    {}_1\phi_{0}\left(\begin{array}{c}
         a\\
         - 
    \end{array};q,z\right)=\frac{(az;q)_{\infty}}{(z;q)_{\infty}}=\sum_{n=0}^{\infty}\frac{(a;q)_{n}}{(q;q)_{n}}z^{n}.
\end{equation}
The $q$-exponential $\e_{q}(z)$ is defined by
\begin{equation*}
    \e_{q}(z)=\sum_{n=0}^{\infty}\frac{z^n}{(q;q)_{n}}={}_1\phi_{0}\left(\begin{array}{c}
         0\\
         - 
    \end{array};q,\minus z\right)=\frac{1}{(z;q)_{\infty}}.
\end{equation*}
Another $q$-analogue of the classical exponential function is
\begin{equation*}
    \E_{q}(z)=\sum_{n=0}^{\infty}q^{\binom{n}{2}}\frac{z^n}{(q;q)_{n}}={}_1\phi_{1}\left(\begin{array}{c}
         0\\
         0 
    \end{array};q,\minus z\right)=(\minus z;q)_{\infty}.
\end{equation*}
The $q$-differential operator $D_{q}$ is defined by:
\begin{equation*}
    D_{q}f(x)=\frac{f(x)-f(qx)}{x}
\end{equation*}
and the Leibniz rule for $D_{q}$
\begin{equation}\label{eqn_leibniz}
    D_{q}^{n}\{f(x)g(x)\}=\sum_{k=0}^{n}q^{k(k\minus n)}\qbinom{n}{k}_{q}D_{q}^{k}\{f(x)\}D_{q}^{n\minus k}\{g(q^{k}x)\}.
\end{equation}
From Eq.(\ref{iden2}), for $k\geq1$,
\begin{equation}
    D_{q}^nx^{\minus k}=(q;q)_{n}\qbinom{-k}{n}_{q}x^{-n-k}=(\minus1)^{n}q^{\minus kn\minus\binom{n}{2}}(q;q)_{n}\qbinom{n+k-1}{n}_{q}x^{\minus k\minus n}.\label{iden9}
\end{equation}
In \cite{orozco}, the author defined, for all $u\in\C$, the $u$-deformed $q$-exponential function as
\begin{equation}\label{def_exp}
    \e_{q}(z,u)=
    \begin{cases}
        \sum_{n=0}^{\infty}u^{\binom{n}{2}}\frac{z^{n}}{(q;q)_{n}}&\text{ if }u\neq0;\\
        1+\frac{z}{1\minus q}&\text{ if }u=0.
    \end{cases}
\end{equation}
Some $u$-deformed $q$-exponential functions are
\begin{align*}
    \e_{q}(z,1)&=e_{q}(z)=\frac{1}{(z;q)_{\infty}},\ \vert z\vert<1,\\
    \e_{q}(\minus z,q)&=\E_{q}(\minus z)=(z;q)_{\infty},\ z\in\C,\\
    \e_{q}(z,\sqrt{q})&=\mathcal{E}_{q}(z)=\sum_{n=0}^{\infty}q^{\frac{1}{2}\binom{n}{2}}\frac{z^n}{(q;q)_{n}}={}_{1}\phi_{1}\left(\begin{array}{c}
         0\\
         -\sqrt{q}
    \end{array};\sqrt{q},\minus z\right),\ z\in\C,\\
    \e_{q}(qz,q^2)&=\mathcal{R}_{q}(z)=\sum_{n=0}^{\infty}q^{n^2}\frac{z^n}{(q;q)_{n}}\\
    &=\frac{(zq^5,z^2q^2,z^2q^3;q^5)_{\infty}}{(zq;q)_{\infty}}{}_{3}\phi_{2}\left(\begin{array}{c}
         z/q,z,zq\\
         z^2q^2,z^2q^3
    \end{array};q^5,zq^5\right),\\
    \e_{q}(z,q^{b+1/2})&=\mathcal{E}_{q,b}(z)=\sum_{n=0}^{\infty}q^{(b+1/2)\binom{n}{2}}\frac{z^n}{(q;q)_{n}},\ b\geq1,\\
    &={}_{0}\phi_{2b}\left(\begin{array}{c}
         -\\
         -\sqrt{q},\mathbf{0}_{2b-1}
    \end{array};\sqrt{q},(-1)^{2b+1}z\right),\ z\in\C,,\\
    \e_{q}(z,q^{b})&=\mathcal{R}_{q,b}(z)=\sum_{n=0}^{\infty}q^{b\binom{n}{2}}\frac{z^n}{(q;q)_{n}},\ b\geq1,\\
    &={}_{0}\phi_{b-1}\left(\begin{array}{c}
         -\\
         \mathbf{0}_{b-1}
    \end{array};q,(-1)^{b-1}z\right),\ z\in\C,
\end{align*}
where $\mathcal{E}_{q}(z)$ is the Exton $q$-exponential function and $\mathcal{R}_{q}(z)$ is the Rogers-Ramanujan function and
\begin{equation}
    \mathbf{0}_{n}=
    \begin{cases}
        (\overbrace{0,\ldots,0)}^n,&\text{ if }n>0;\\
        \emptyset,&\text{ if }n=0.
    \end{cases}
\end{equation}
Next, we present the following lemmas that will help us derive the $q$-identities in this paper.
\begin{lemma}\label{lemma_prop_exp}
For $0<\vert q\vert<1$ and $x\neq0$
\begin{align*}
    \frac{1}{(q^nx;q)_{\infty}}&=\frac{(x;q)_{n}}{(x;q)_{\infty}},\\
    \frac{1}{(q^{-n}x;q)_{\infty}}&=\frac{(\minus1)^nq^{\binom{n+1}{2}}x^{-n}}{(q/x;q)_{n}(x;q)_{\infty}}.
\end{align*}    
\end{lemma}

\begin{lemma}\label{lemma_prop_exp2}
For $0<\vert q\vert<1$ and $x\neq0$
\begin{align*}
    (q^{n}x,q)_{\infty}&=\frac{(x,q)_{\infty}}{(x;q)_{n}},\\
    (q^{\minus n}x;q)_{\infty}&=(-1)^nx^nq^{-\binom{n+1}{2}}(q/x;q)_{n}(x;q)_{\infty}.
\end{align*}    
\end{lemma}

\begin{lemma}\label{lemma_prop_exp4}
For $0<\vert q\vert<1$ and $z\neq0$
\begin{align*}
\mathcal{R}_{q}(q^{5n}z)&=\frac{(zq^5,z^2q^2,z^2q^3;q^5)_{\infty}(zq;q)_{5n}}{(zq;q)_{\infty}(zq^{5};q^5)_{n}(z^2q^2;q^5)_{2n}(z^2q^3;q^5)_{2n}}\\
&\hspace{2cm}\times{}_{3}\phi_{2}\left(\begin{array}{c}
         zq^{5n-1},zq^{5n},zq^{5n+1}\\
         z^2q^{10n+2},z^2q^{10n+3}
    \end{array};q^5,zq^{5(n+1)}\right).\\   
\mathcal{R}_{q}(q^{-5n}z)&=\frac{(zq^5,z^2q^2,z^2q^3;q^5)_{\infty}}{(zq;q)_{\infty}}q^{5n}z^{4n}\frac{(z^{-1};q)_{n}(z^{-2}q^3;q^5)_{2n}(z^{-2}q^2;q^5)_{2n}}{(z^{-1};q)_{5n}}\\
&\hspace{2cm}\times{}_{3}\phi_{2}\left(\begin{array}{c}
         zq^{-5n-1},zq^{-5n},zq^{-5n+1}\\
         z^2q^{-10n+2},z^2q^{-10n+3}
    \end{array};q^5,zq^{-5(n-1)}\right).
\end{align*}
\end{lemma}

\section{$u$-deformed basic hypergeometric series}

In \cite{orozco} was defined the $u$-$u$-deformed basic hypergeometric series ${}_{r}\Phi_{s}$ as
   \begin{align}
        &{}_{r}\Phi_{s}\left(
    \begin{array}{c}
         a_{1},a_{2},\ldots,a_{r} \\
         b_{1},\ldots,b_{s}
    \end{array}
    ;q,u,z
    \right)\nonumber\\
    &\hspace{3cm}=\sum_{n=0}^{\infty}u^{\binom{n}{2}}\frac{(a_{1},a_{2},\ldots,a_{r};q)_{n}}{(q,b_{1},b_{2},\ldots,b_{s};q)_{n}}\bigg[(\minus1)^nq^{\binom{n}{2}}\bigg]^{1+s-r}z^n,
    \end{align}
where $0<\vert q\vert<1$ and $u\in\C$. Some convergence conditions for $_{r}\Phi_{s}$-series are
\begin{itemize}
    \item $1+s-r\neq0$. If $0<\vert uq^{1+s-r}\vert<1$, then ${}_{r}\Phi_{s}$ is an entire function. If $\vert uq^{1+s-r}\vert=1$, then ${}_{r}\Phi_{s}$ converges for $\vert z\vert<1$. If $\vert uq^{1+s-r}\vert>1$, then ${}_{r}\Phi_{s}$ is divergent.
    \item $1+s-r=0$. If $0<\vert u\vert<1$, then ${}_{r}\Phi_{s}$ is an entire function. If $\vert u\vert=1$, then ${}_{r}\Phi_{s}$ converges for $\vert z\vert<1$. If $\vert u\vert>1$, then ${}_{r}\Phi_{s}$ is divergent.
\end{itemize}
By letting that $u$ take values $q^b$ or $q^{b+1/2}$, we obtain the following basic hypergeometric series:
\begin{itemize}
    \item If $u=q^b$, $b\geq0$, 
\begin{equation}
    {}_{1+r}\Phi_{r}\left(
    \begin{array}{c}
         a_{1},a_{2},\ldots,a_{r+1}\\
         b_{1},\ldots,b_{r}
    \end{array}
    ;q,q^b,z
    \right)={}_{1+r}\phi_{r+b}\left(
    \begin{array}{c}
         a_{1},a_{2},\ldots,a_{r+1} \\
         b_{1},b_{2},\ldots,b_{r}.\mathbf{0}_{b}
    \end{array}
    ;q,(\minus1)^bz
    \right)
\end{equation}
for all $z\in\C$.
\item If $u=q^{b+1/2}$, $b\geq0$
\begin{align}
    &{}_{r+1}\Phi_{r}\left(
    \begin{array}{c}
         a_{1},\ldots,a_{r+1} \\
         b_{1},\ldots,b_{r}
    \end{array}
    ;q,q^{b+1/2},z
    \right)\nonumber\\
    &\hspace{1cm}={}_{2r+2}\phi_{2(r+b)+2}\left(
    \begin{array}{ccc}
         \sqrt{a_{1}},-\sqrt{a_{1}}&,\ldots,&\sqrt{a_{r+1}},-\sqrt{a_{r+1}}\hspace{0.6cm}\\
         \sqrt{b_{1}},-\sqrt{b_{1}}&,\ldots,&\sqrt{b_{r}},-\sqrt{b_{r}},-\sqrt{q},\mathbf{0}_{2b+1}
    \end{array}
    ;\sqrt{q},\minus z
    \right)
\end{align}
for all $z\in\C$. 
\end{itemize}

\section{$u$-deformed homogeneous functions}

\subsection{Definition}

\begin{definition}\label{def_function}
For all $\alpha\in\C$ and $0<q<1$, we define the $u$-deformed homogeneous functions as
\begin{equation}\label{eqn_nbs}
    \rr_{\alpha}(x,y;u|q)=\sum_{n=0}^{\infty}\qbinom{\alpha}{n}_{q}u^{\binom{n}{2}}x^{\alpha\minus n}y^{n}.
\end{equation}
\end{definition}
From Eq.(\ref{iden1}) we have the representation in $q$-series  for the function $\rr_{\alpha}(x,y;u|q)$.
\begin{equation}\label{qseries_rr}
    \rr_{\alpha}(x,y;u|q)=x^\alpha\sum_{n=0}^{\infty}(u/q)^{\binom{n}{2}}\frac{(q^{-\alpha};q)_{n}}{(q;q)_{n}}(yq^\alpha/x)^{n}.
\end{equation}
If $\alpha=-n$, then by the Eq.(\ref{iden2}) the function $\rr_{\minus n}(x,y;u|q)$ is the generating function of the $q$-binomial coefficients $\qbinom{n+k-1}{k}_{q}$
\begin{equation}\label{rr_as_GF}
    \rr_{\minus n}(x,y;u|q)=x^{-n}\sum_{k=0}^{\infty}\qbinom{n+k-1}{k}_{q}(u/q)^{\binom{k}{2}}(-y/q^nx)^k.
\end{equation}
\begin{theorem}\label{theo_conver}
Take $x,y,\alpha\in\C$, such that $x\neq0$ and $\alpha$ is not a non-negative integer. The function $\rr_{\alpha}(x,y;u|q)$ is an entire function in the variable $y$ if $0<\vert u\vert<\vert q\vert<1$. If $\vert u\vert=\vert q\vert<1$, then $\rr_{\alpha}(x,y;u|q)$ converges for $\vert y/x\vert<\vert q^{-\alpha}\vert$. Finally, $\rr_{\alpha}(x,y;u|q)$ converges in $y=0$ if $\vert u\vert>\vert q\vert$.
\end{theorem}
\begin{proof}
From ratio's test applied to Eq.(\ref{qseries_rr})
\begin{align}\label{test_ratio}
    \bigg\vert\frac{q^\alpha y}{x}\bigg\vert\lim_{n\rightarrow\infty}\bigg\vert\left(\frac{u}{q}\right)^n\frac{1-q^{-\alpha+n}}{1-q^{n+1}}\bigg\vert&=\lim_{n\rightarrow\infty}\bigg\vert\left(\frac{u}{q}\right)^n\bigg\vert=
    \begin{cases}
        0,&\text{ if }\vert u\vert<\vert q\vert;\\
        \big\vert\frac{q^\alpha y}{x}\big\vert,&\text{ if }\vert u\vert=\vert q\vert;\\
        \infty,&\text{ if }\vert u\vert>\vert q\vert.
    \end{cases}
\end{align}
The statements of the theorem are deduced from Eq.(\ref{test_ratio}).
\end{proof}

When $\alpha=n\in\N$, we obtain the $u$-deformed homogeneous polynomials defined in \cite{orozco}
\begin{equation*}
    \rr_{n}(x,y;u|q)=\sum_{k=0}^{n}\qbinom{n}{k}_{q}u^{\binom{k}{2}}x^{n\minus k}y^{k}.
\end{equation*}
From Definition \ref{def_function} and Theorem \ref{theo_conver}, we obtain the following specializations.
\begin{itemize}
\item If $u=q$, we define the Cauchy function as
\begin{align}\label{cauchy_alpha}
    \P_{\alpha}(x,y|q)&=\rr_{\alpha}(x,-y;q|q)=\sum_{k=0}^{\infty}\qbinom{\alpha}{k}_{q}(\minus1)^kq^{\binom{k}{2}}x^{\alpha\minus k}y^{k}
\end{align}
which is convergent for all $\vert y/x\vert<\vert q^{-\alpha}\vert$ if $0<\vert q\vert<1$. A representation in infinite product for the function in Eq.(\ref{cauchy_alpha}) is
\begin{equation}
    \P_{\alpha}(x,y|q)=\frac{\P_{\infty}(x,y|q)}{\P_{\infty}(x,q^\alpha y|q)},
\end{equation}
where 
\begin{equation}
    \P_{\infty}(x,y|q)=\prod_{k=0}^{\infty}(x-q^ky).
\end{equation}
From Eq.(\ref{rr_as_GF}), 
\begin{equation}
    \P_{\minus n}(x,y|q)=x^{-n}\sum_{k=0}^{\infty}\qbinom{n+k-1}{k}_{q}(y/q^nx)^k.
\end{equation}

\item If $u=q^2$, we define the Stieltjes-Wigert function as
\begin{align}
    \s_{\alpha}(x,y;q)&=\rr_{\alpha}(x,qy;q^2|q)\nonumber\\
    &=\sum_{k=0}^{\infty}\qbinom{\alpha}{k}_{q}(\minus1)^kq^{k^2}x^{\alpha\minus k}y^{k},
\end{align}
which is an entire function in the variable $y$.
From Eq.(\ref{rr_as_GF})
\begin{equation}
    \s_{\minus n}(x,y;q)=x^{-n}\sum_{k=0}^{\infty}\qbinom{n+k-1}{k}_{q}q^{\binom{k}{2}}(-y/q^{n-1}x)^k.
\end{equation}
In general, if $u=q^{b}$, $b\geq1$, we define the Stieltjes-Wigert function of order $b$ as
\begin{align}
    \s_{\alpha,b}(x,y;q)&=\rr_{\alpha}(x,y;q^b|q)\nonumber\\
    &=\sum_{k=0}^{\infty}\qbinom{\alpha}{k}_{q}(\minus1)^kq^{b\binom{k}{2}}x^{\alpha\minus k}y^{k},
\end{align}
which is an entire function in the variable $y$.
From Eq.(\ref{rr_as_GF})
\begin{equation}
    \s_{\minus n,b}(x,y;q)=x^{-n}\sum_{k=0}^{\infty}\qbinom{n+k-1}{k}_{q}q^{(b-1)\binom{k}{2}}(-y/q^{n-1}x)^k.
\end{equation}

\item If $u=q^{b+1/2}$, $b\geq1$, we define the Exton function of order $b$ as
\begin{align}
    \E_{\alpha,b}(x,y;q)&=\rr_{\alpha}(x,y;q^{b+1/2}|q)=\sum_{k=0}^{\infty}\qbinom{\alpha}{k}_{q}q^{(b+\frac{1}{2})\binom{k}{2}}x^{\alpha\minus k}y^{k}\nonumber\\
    &=\sum_{k=0}^{\infty}\qbinom{\alpha}{k}_{\sqrt{q}}\frac{(\minus\sqrt{q};\sqrt{q})_{\alpha}q^{(b+\frac{1}{2})\binom{k}{2}}x^{\alpha\minus k}y^{k}}{(\minus\sqrt{q};\sqrt{q})_{k}(-\sqrt{q};\sqrt{q})_{\alpha\minus k}},
\end{align}
which is an entire function in the variable $y$. From Eq.(\ref{rr_as_GF})
\begin{equation}
    \E_{\minus n,b}(x,y;q)=x^{-n}\sum_{k=0}^{\infty}\qbinom{n+k-1}{k}_{q}q^{(b-\frac{1}{2})\binom{k}{2}}(-y/q^{n-1}x)^k.
\end{equation}
\end{itemize}

\subsection{Basic properties}
Similarly to what was done in \cite{orozco}, we obtain some basic properties for the $u$-deformed homogeneous functions.
\begin{itemize}
\item 
Representation with $u$-deformed basic hypergeometric series. For all $\alpha\in\C$ and for $x\neq0$,
    \begin{equation}\label{eqn_rep_bhs}
    \rr_{\alpha}(x,y;u|q)=x^\alpha{}_{2}\Phi_{0}\left(
    \begin{array}{c}
         q^{\minus\alpha},0 \\
         -
    \end{array}
    ;q,u,q^\alpha y/x    
    \right).
\end{equation}
\item 
Recurrence identities,
\begin{align*}
    \rr_{\alpha+1}(x,y;u|q)&=x\rr_{\alpha}(x,qy;u|q)+y\rr_{\alpha}(x,uy;u|q)
\end{align*}
and
\begin{equation*}
    \rr_{\alpha+1}(x,y;u|q)=x\rr_{\alpha}(x,y;u|q)+y\rr_{\alpha}(qx,uy;u|q).    
\end{equation*}
\item 
$q$-derivatives,
\begin{align*}
    \mathbf{D}_{q}\{\rr_{\alpha}(x,a;u|q)\}&=(1- q^\alpha)\rr_{\alpha\minus1}(x,a;u|q),
\end{align*}
\begin{align*}
    \mathbf{D}_{q}\{\rr_{\alpha}(a,x;u|q)\}&=(1-q^\alpha)\rr_{\alpha\minus1}(a,ux;u|q),
\end{align*}
and
\begin{equation*}
    \mathbf{D}_{q}\{\rr_{\alpha}(a,-x;u|q)\}=-(1- q^\alpha)\rr_{\alpha\minus1}(a,-ux;u|q).    
\end{equation*}
For $k\geq2$,
\begin{align}
    \mathbf{D}^{k}_{q}\{\rr_{\alpha}(x,a;u|q)\}&=(q;q)_{k}\qbinom{\alpha}{k}_{q}\rr_{\alpha\minus k}(x,a;u|q),\label{eqn_difk_rsf1}
\end{align}
and
\begin{align}
    \mathbf{D}^{k}_{q}\{\rr_{\alpha}(a,x;u|q)\}&=u^{\binom{k}{2}}(q;q)_{k}\qbinom{\alpha}{k}_{q}\rr_{\alpha\minus k}(a,u^kx;u|q).\label{eqn_difk_rsf2}.
\end{align}
\end{itemize}
In \cite{orozco}, for all $u\in\C$, was define the $u$-deformed $q$-exponential operator $\E(yD_{q}|u)$ by letting 
\begin{equation*}
    \E(yD_{q}|u)=\sum_{n=0}^{\infty}u^{\binom{n}{2}}\frac{(yD_{q})^{n}}{(q;q)_{n}}. 
\end{equation*}
This operator is a generalization of the $q$-exponential operators of Chen \cite{chen1} and Saad \cite{saad}. The $u$-deformed homogeneous function $\rr_{\alpha}(x,y;u|q)$ can be represented by the $u$-deformed $q$-exponential operator $\E(yD_{q}|u)$ as follows.
\begin{proposition}\label{prop_translation}
For all $\alpha,u\in\C$
    \begin{align*}
        \E(yD_{q}|u)\big\{x^{\alpha}\big\}&=\rr_{\alpha}(x,y;u|q).
    \end{align*}
\end{proposition}

\section{Generating functions for the function $\rr_{\minus n}(x,y;u|q)$}

\subsection{Some $q$-series}

\begin{lemma}\label{lemma_nder_basic}
For all $n\geq0$ and $x\neq0$,
\begin{align}
        D_{q}^{n}\{\e_{q}(a/x,u)\}&=(\minus1)^nq^{\minus\binom{n}{2}}x^{\minus n}\sum_{i=0}^{n}\qbinom{n}{i}_{q}(\minus1)^iq^{\binom{i}{2}}\e_{q}(aq^{i\minus n}/x,u).
    \end{align}
\end{lemma}
\begin{proof}
Follows from Eqs. (\ref{iden9}) and  (\ref{iden1})
    \begin{align*}
        D_{q}^n\{\e_{q}(a/x,u)\}&=\sum_{k=0}^{\infty}u^{\binom{k}{2}}\frac{a^k}{(q;q)_{k}}D_{q}^{n}\{x^{\minus k}\}\\
        &=(\minus1)^nq^{\minus\binom{n}{2}}x^{\minus n}\sum_{k=0}^{\infty}u^{\binom{k}{2}}\frac{(q^{k};q)_{n}}{(q;q)_{k}}(q^{\minus n}a/x)^{k}\\
        &=(\minus1)^nq^{\minus\binom{n}{2}}x^{\minus n}\sum_{k=0}^{\infty}u^{\binom{k}{2}}\frac{(q^{\minus n}a/x)^{k}}{(q;q)_{k}}\sum_{i=0}^{n}\qbinom{n}{i}_{q}(\minus1)^{i}q^{ki+\binom{i}{2}}\\
        &=(\minus1)^nq^{\minus\binom{n}{2}}x^{\minus n}\sum_{i=0}^{n}\qbinom{n}{i}_{q}(\minus1)^iq^{\binom{i}{2}}\sum_{k=0}^{\infty}u^{\binom{k}{2}}\frac{(q^{i\minus n}a/x)^{k}}{(q;q)_{k}}\\
        &=(\minus1)^nq^{\minus\binom{n}{2}}x^{\minus n}\sum_{i=0}^{n}\qbinom{n}{i}_{q}(\minus1)^iq^{\binom{i}{2}}\e_{q}(aq^{i\minus n}/x,u).
        \end{align*}
The claim is reached.        
\end{proof}

\begin{theorem}\label{theo_qoper_exp}
If $0<\vert u\vert\leq q$ and $x\neq0$, then
    \begin{equation}
        \E(yD_{q}|u)\{\e_{q}(a/x,v)\}=\sum_{i=0}^{\infty}u^{\binom{i}{2}}\frac{(y/x)^{i}}{(q;q)_{i}}\sum_{n=0}^{\infty}(q^{\minus1}u)^{\binom{n}{2}}\frac{(-u^iy/x)^{n}}{(q;q)_{n}}\e_{q}(aq^{\minus n}/x,v).
    \end{equation}
\end{theorem}
\begin{proof}
From Lemma \ref{lemma_nder_basic}, we have that
    \begin{align*}
        \E(yD_{q}|u)\{\e_{q}(a/x,v)\}&=\sum_{n=0}^{\infty}u^{\binom{n}{2}}\frac{y^n}{(q;q)_{n}}D_{q}^{n}\{\e_{q}(a/x,v)\}\\
        &=\sum_{n=0}^{\infty}(q^{\minus1}u)^{\binom{n}{2}}\frac{(\minus y/x)^n}{(q;q)_{n}}\sum_{i=0}^{n}\qbinom{n}{i}_{q}(\minus1)^iq^{\binom{i}{2}}\e_{q}(aq^{i\minus n}/x,v)\\
        &=\sum_{i=0}^{\infty}u^{\binom{i}{2}}\frac{(y/x)^{i}}{(q;q)_{i}}\sum_{n=0}^{\infty}(q^{\minus1}u)^{\binom{n}{2}}\frac{(-u^iy/x)^{n}}{(q;q)_{n}}\e_{q}(aq^{\minus n}/x,v),
    \end{align*}
such as is claimed.    
\end{proof}

\begin{theorem}\label{theo_gdf}
If $0<\vert u\vert\leq q$ and $x\neq0$, then
    \begin{multline}
        \sum_{n=0}^{\infty}v^{\binom{n}{2}}\rr_{\minus n}(x,y;u|q)\frac{z^n}{(q;q)_{n}}\\
        =\sum_{i=0}^{\infty}u^{\binom{i}{2}}\frac{(y/x)^{i}}{(q;q)_{i}}\sum_{n=0}^{\infty}(u/q)^{\binom{n}{2}}\frac{(-u^iy/x)^{n}}{(q;q)_{n}}\e_{q}(zq^{\minus n}/x,v).
    \end{multline}
\end{theorem}
\begin{proof}
The result follows from Theorem \ref{theo_qoper_exp} by noting that 
    \begin{align*}
        &\sum_{n=0}^{\infty}v^{\binom{n}{2}}\rr_{\minus n}(x,y;u|q)\frac{z^n}{(q;q)_{n}}\\
        &\hspace{1cm}=\E(yD_{q}|u)\left\{\sum_{n=0}^{\infty}v^{\binom{n}{2}}\frac{(z/x)^n}{(q;q)_{n}}\right\}\\
        &\hspace{1cm}=\E(yD_{q})\left\{\e_{q}(z/x,v)\right\}\\
        &\hspace{1cm}=\sum_{i=0}^{\infty}u^{\binom{i}{2}}\frac{(x^{\minus1}y)^{i}}{(q;q)_{i}}\sum_{n=0}^{\infty}(q^{\minus1}u)^{\binom{n}{2}}\frac{(-u^ix^{\minus1}y)^{n}}{(q;q)_{n}}\e_{q}(zq^{\minus n}/x,v).
    \end{align*}
\end{proof}

\begin{theorem}\label{theo_operE1}
If $0<\vert u\vert<1$ and $x\neq0$, $a\neq0$, then
    \begin{equation}
        \E(yD_{q}|u)\left\{\frac{1}{(a/x;q)_{\infty}}\right\}=\frac{1}{(a/x;q)_{\infty}}\sum_{i=0}^{\infty}u^{\binom{i}{2}}\frac{(y/x)^{i}}{(q;q)_{i}}{}_{2}\Phi_{1}\left(
        \begin{array}{c}
             0,0 \\
             qx/a
        \end{array};
        q,u,u^iqy/a
        \right).
    \end{equation}
\end{theorem}
\begin{proof}
Use $v=1$ in Theorem \ref{theo_qoper_exp} and then Lemma \ref{lemma_prop_exp}. Hence
    \begin{align*}
        &\E(yD_{q}|u)\left\{\frac{1}{(a/x;q)_{\infty}}\right\}\\
        &=\sum_{i=0}^{\infty}u^{\binom{i}{2}}\frac{(y/x)^{i}}{(q;q)_{i}}\sum_{n=0}^{\infty}(q^{\minus1}u)^{\binom{n}{2}}\frac{(-u^iy/x)^{n}}{(q;q)_{n}}\frac{1}{(aq^{\minus n}/x;q)_{\infty}}\\
        &=\sum_{i=0}^{\infty}u^{\binom{i}{2}}\frac{(y/x)^{i}}{(q;q)_{i}}\sum_{n=0}^{\infty}(q^{\minus1}u)^{\binom{n}{2}}\frac{(-u^iy/x)^{n}}{(q;q)_{n}}\frac{(-1)^nq^{\binom{n+1}{2}}(a/x)^{-n}}{(qx/a;q)_{n}(a/x;q)_{\infty}}\\
        &=\frac{1}{(a/x;q)_{\infty}}\sum_{i=0}^{\infty}u^{\binom{i}{2}}\frac{(y/x)^{i}}{(q;q)_{i}}\sum_{n=0}^{\infty}u^{\binom{n}{2}}\frac{(u^iqy/a)^{n}}{(qx/a;q)_{n}(q;q)_{n}}\\
        &=\frac{1}{(a/x;q)_{\infty}}\sum_{i=0}^{\infty}u^{\binom{i}{2}}\frac{(y/x)^{i}}{(q;q)_{i}}{}_{2}\Phi_{1}\left(
        \begin{array}{c}
             0,0 \\
             qx/a
        \end{array};
        q,u,u^iqy/a
        \right).
    \end{align*}
The claim is proved.    
\end{proof}

\begin{theorem}\label{theo_gdfE1}
If $0<\vert u\vert<1$ and $x\neq0$, $z\neq0$, then
    \begin{align}
        \sum_{n=0}^{\infty}\rr_{\minus n}(x,y;u|q)\frac{z^n}{(q;q)_{n}}=\frac{1}{(z/x;q)_{\infty}}\sum_{i=0}^{\infty}u^{\binom{i}{2}}\frac{(y/x)^{i}}{(q;q)_{i}}{}_{2}\Phi_{1}\left(
        \begin{array}{c}
             0,0 \\
             qx/z
        \end{array};
        q,u,u^iqy/z
        \right).
    \end{align}
\end{theorem}
\begin{proof}
Followed by applying Theorem \ref{theo_gdf} with $v=1$ and then Theorem \ref{theo_operE1}.    
\end{proof}

\begin{theorem}\label{theo_operE2}
If $0<\vert u\vert\leq q^2$ and $x\neq0$, $a\neq0$, then
    \begin{multline}
        \E(yD_{q}|u)\left\{(a/x;q)_{\infty}\right\}\\
        =(a/x;q)_{\infty}\sum_{i=0}^{\infty}u^{\binom{i}{2}}\frac{(y/x)^{i}}{(q;q)_{i}}{}_{1}\Phi_{0}\left(
        \begin{array}{c}
             qx/a \\
             -
        \end{array};
        q,u/q^2,au^iy/qx^2
        \right).
    \end{multline}
\end{theorem}
\begin{proof}
Use $v=q$ and replace $a$ by $-a$ in Theorem \ref{theo_qoper_exp} and finally use the Lemma \ref{lemma_prop_exp2}. Hence
    \begin{align*}
        &\E(yD_{q})\left\{(a/x;q)_{\infty}\right\}\\
        &=\sum_{i=0}^{\infty}u^{\binom{i}{2}}\frac{(y/x)^{i}}{(q;q)_{i}}\sum_{n=0}^{\infty}(u/q)^{\binom{n}{2}}\frac{(-u^iy/x)^{n}}{(q;q)_{n}}(aq^{\minus n}/x;q)_{\infty}\\
        &=\sum_{i=0}^{\infty}u^{\binom{i}{2}}\frac{(y/x)^{i}}{(q;q)_{i}}\sum_{n=0}^{\infty}(u/q)^{\binom{n}{2}}\frac{(-u^iy/x)^{n}}{(q;q)_{n}}(-1)^n(a/x)^nq^{-\binom{n+1}{2}}(qx/a;q)_{n}(a/x;q)_{\infty}\\
        &=(a/x;q)_{\infty}\sum_{i=0}^{\infty}u^{\binom{i}{2}}\frac{(y/x)^{i}}{(q;q)_{i}}\sum_{n=0}^{\infty}(u/q^2)^{\binom{n}{2}}\frac{(qx/a;q)_{n}(au^iy/qx^2)^{n}}{(q;q)_{n}}\\
        &=(a/x;q)_{\infty}\sum_{i=0}^{\infty}u^{\binom{i}{2}}\frac{(x^{\minus1}y)^{i}}{(q;q)_{i}}{}_{1}\Phi_{0}\left(
        \begin{array}{c}
             qx/a \\
             -
        \end{array};
        q,u/q^2,au^iy/qx^2
        \right).
    \end{align*}
The claim is proved.    
\end{proof}

\begin{theorem}\label{theo_gdfE2}
If $0<\vert u\vert\leq q^2$ and $x\neq0$, $z\neq0$, then
    \begin{multline}
        \sum_{n=0}^{\infty}(-1)^nq^{\binom{n}{2}}\rr_{\minus n}(x,y;u|q)\frac{z^n}{(q;q)_{n}}\\
        =(z/x;q)_{\infty}\sum_{i=0}^{\infty}u^{\binom{i}{2}}\frac{(y/x)^{i}}{(q;q)_{i}}{}_{1}\Phi_{0}\left(
        \begin{array}{c}
             qx/z \\
             -
        \end{array};
        q,u/q^2,zu^iy/qx^2
        \right).
    \end{multline}
\end{theorem}
\begin{proof}
Follow from Theorem \ref{theo_gdf} with $v=q$ and from Theorem \ref{theo_operE2}.    
\end{proof}

\subsection{Some $q^{b+1/2}$-series}

\begin{lemma}\label{lemma_Ebnder}
For all $n\geq0$ and $x\neq0$,
\begin{multline}
        D_{q^{b+1/2}}^{n}\{\mathcal{E}_{q,b}(a/x)\}\\
        =(\minus1)^nq^{\minus(b+1/2)\binom{n}{2}}x^{\minus n}\sum_{i=0}^{n}\qbinom{n}{i}_{q^{b+1/2}}(\minus1)^iq^{(b+1/2)\binom{i}{2}}\mathcal{E}_{q,b}(aq^{(b+1/2)(i-n)}/x).
    \end{multline}
\end{lemma}
\begin{proof}
Replace $q$ by $q^{b+1/2}$ and set $u=q^{b+1/2}$ in Lemma \ref{lemma_nder_basic}.    
\end{proof}

\begin{theorem}\label{theo_Eoper}
If $0<\vert u\vert\leq\sqrt{q}$ and $x\neq0$, then
    \begin{multline}
        \E(yD_{\sqrt{q}}|u)\left\{\mathcal{E}_{q}(a/x)\right\}\\
        =\sum_{i=0}^{\infty}u^{\binom{i}{2}}\frac{(y/x)^{i}}{(\sqrt{q};\sqrt{q})_{i}}\sum_{n=0}^{\infty}(q^{-\frac{1}{2}}u)^{\binom{n}{2}}\frac{(-q^{-\frac{i}{2}}u^iy/x)^{n}}{(\sqrt{q};\sqrt{q})_{n}}{}_{1}\phi_{1}\left(\begin{array}{c}
         0\\
         -\sqrt{q}
    \end{array};\sqrt{q},\minus q^{-n}a/x\right).
    \end{multline}
\end{theorem}
\begin{proof}
Use Lemma \ref{lemma_Ebnder} with $b=0$ and then the representation of the function $\mathcal{E}_{q}(z)$. 
    \begin{align*}
        &\E(yD_{\sqrt{q}}|u)\left\{\mathcal{E}_{q}(a/x)\right\}\\
        &=\sum_{n=0}^{\infty}u^{\binom{n}{2}}\frac{y^n}{(\sqrt{q};\sqrt{q})_{n}}D_{\sqrt{q}}^{n}\{\mathcal{E}_{q}(a/x)\}\\
        &=\sum_{n=0}^{\infty}(q^{-\frac{1}{2}}u)^{\binom{n}{2}}\frac{(\minus y/x)^n}{(\sqrt{q};\sqrt{q})_{n}}\sum_{i=0}^{n}\qbinom{n}{i}_{\sqrt{q}}(\minus1)^iq^{\frac{1}{2}\binom{i}{2}}\mathcal{E}_{q}(aq^{\frac{i\minus n}{2}}/x)\\
        &=\sum_{i=0}^{\infty}u^{\binom{i}{2}}\frac{(y/x)^{i}}{(\sqrt{q};\sqrt{q})_{i}}\sum_{n=0}^{\infty}(q^{-\frac{1}{2}}u)^{\binom{n}{2}}\frac{(-q^{-\frac{i}{2}}u^iy/x)^{n}}{(\sqrt{q};\sqrt{q})_{n}}\mathcal{E}_{q}(aq^{\minus n}/x)\\
        &=\sum_{i=0}^{\infty}u^{\binom{i}{2}}\frac{(y/x)^{i}}{(\sqrt{q};\sqrt{q})_{i}}\sum_{n=0}^{\infty}(q^{-\frac{1}{2}}u)^{\binom{n}{2}}\frac{(-q^{-\frac{i}{2}}u^iy/x)^{n}}{(\sqrt{q};\sqrt{q})_{n}}{}_{1}\phi_{1}\left(\begin{array}{c}
         0\\
         -\sqrt{q}
    \end{array};\sqrt{q},\minus q^{-n}a/x\right).
    \end{align*}
The claim is proved.    
\end{proof}

\begin{theorem}\label{theo_Eboper}
If $0<\vert u\vert\leq q^{b+1/2}$, $b\geq1$ and $x\neq0$, then
    \begin{multline}
        \E(yD_{q^{b+1/2}}|u)\left\{\mathcal{E}_{q,b}(a/x)\right\}\\
        =\sum_{i=0}^{\infty}u^{\binom{i}{2}}\frac{(y/x)^{i}}{(\sqrt{q}q^b;\sqrt{q}q^b)_{i}}\sum_{n=0}^{\infty}(q^{-(b+\frac{1}{2})}u)^{\binom{n}{2}}\frac{(-q^{-(b+\frac{1}{2})i}u^iy/x)^{n}}{(\sqrt{q}q^b;\sqrt{q}q^b)_{n}}\\
        \times{}_{0}\phi_{2b}\left(\begin{array}{c}
         -\\
         -\sqrt{q},\mathbf{0}_{2b-1}
    \end{array};\sqrt{q},(-q^{n/2})^{2b+1}a/x\right).
    \end{multline}
\end{theorem}
\begin{proof}
Use Lemma \ref{lemma_Ebnder} with $b\geq1$ and then the representation of the function $\mathcal{E}_{q,b}(z)$.
    \begin{align*}
        &\E(yD_{q^{b+\frac{1}{2}}}|u)\left\{\mathcal{E}_{q,b}(a/x)\right\}\\
        &=\sum_{n=0}^{\infty}u^{\binom{n}{2}}\frac{y^n}{(\sqrt{q}q^b;\sqrt{q}q^b)_{n}}D_{q^{b+\frac{1}{2}}}^{n}\{\mathcal{E}_{q,b}(a/x)\}\\
        &=\sum_{n=0}^{\infty}(q^{-(b+\frac{1}{2})}u)^{\binom{n}{2}}\frac{(-y/x)^n}{(\sqrt{q}q^b;\sqrt{q}q^b)_{n}}\sum_{i=0}^{n}\qbinom{n}{i}_{q^{b+\frac{1}{2}}}(\minus1)^iq^{(b+\frac{1}{2})\binom{i}{2}}\mathcal{E}_{q,b}(aq^{(b+\frac{1}{2})(i-n)}/x)\\
        &=\sum_{i=0}^{\infty}u^{\binom{i}{2}}\frac{(y/x)^{i}}{(\sqrt{q}q^b;\sqrt{q}q^b)_{i}}\sum_{n=0}^{\infty}(q^{-(b+\frac{1}{2})}u)^{\binom{n}{2}}\frac{(-q^{-(b+\frac{1}{2})i}u^iy/x)^{n}}{(\sqrt{q}q^b;\sqrt{q}q^b)_{n}}\mathcal{E}_{q,b}(aq^{\minus(b+\frac{1}{2})n}/x).\\
        &=\sum_{i=0}^{\infty}u^{\binom{i}{2}}\frac{(y/x)^{i}}{(\sqrt{q}q^b;\sqrt{q}q^b)_{i}}\sum_{n=0}^{\infty}(q^{-(b+\frac{1}{2})}u)^{\binom{n}{2}}\frac{(-q^{-(b+\frac{1}{2})i}u^iy/x)^{n}}{(\sqrt{q}q^b;\sqrt{q}q^b)_{n}}\\
        &\hspace{5cm}\times{}_{0}\phi_{2b}\left(\begin{array}{c}
         -\\
         -\sqrt{q},\mathbf{0}_{2b-1}
    \end{array};\sqrt{q},(-q^{n/2})^{2b+1}a/x\right).
    \end{align*}
The claim is proved.    
\end{proof}

\begin{theorem}\label{theo_Egdf}
If $0<\vert u\vert\leq\sqrt{q}$ and $x\neq0$, then
    \begin{multline}
        \sum_{n=0}^{\infty}q^{\frac{1}{2}\binom{n}{2}}\rr_{\minus n}(x,y;u|\sqrt{q})\frac{z^n}{(q;q)_{n}}\\
        \sum_{i=0}^{\infty}u^{\binom{i}{2}}\frac{(y/x)^{i}}{(\sqrt{q};\sqrt{q})_{i}}\sum_{n=0}^{\infty}(q^{-\frac{1}{2}}u)^{\binom{n}{2}}\frac{(-q^{-\frac{i}{2}}u^iy/x)^{n}}{(\sqrt{q};\sqrt{q})_{n}}{}_{1}\phi_{1}\left(\begin{array}{c}
         0\\
         -\sqrt{q}
    \end{array};\sqrt{q},\minus q^{-n}z/x\right).
    \end{multline}
\end{theorem}
\begin{proof}
Follow from Theorem \ref{theo_Eoper}. Hence
    \begin{align*}
        &\sum_{n=0}^{\infty}q^{\frac{1}{2}\binom{n}{2}}\rr_{\minus n}(x,y;u|\sqrt{q})\frac{z^n}{(q;q)_{n}}\\
        &\hspace{1cm}=\E(yD_{\sqrt{q}}|u)\left\{\sum_{n=0}^{\infty}q^{\frac{1}{2}\binom{n}{2}}\frac{(z/x)^n}{(q;q)_{n}}\right\}\\
        &\hspace{1cm}=\E(yD_{\sqrt{q}}|u)\left\{\mathcal{E}_{q}(z/x)\right\}\\
        &\hspace{1cm}=\sum_{i=0}^{\infty}u^{\binom{i}{2}}\frac{(y/x)^{i}}{(\sqrt{q};\sqrt{q})_{i}}\sum_{n=0}^{\infty}(q^{-\frac{1}{2}}u)^{\binom{n}{2}}\frac{(-q^{-\frac{i}{2}}u^iy/x)^{n}}{(\sqrt{q};\sqrt{q})_{n}}{}_{1}\phi_{1}\left(\begin{array}{c}
         0\\
         -\sqrt{q}
    \end{array};\sqrt{q},\minus q^{-n}z/x\right).
    \end{align*}
The claim is reached.    
\end{proof}

\begin{theorem}\label{theo_Ebgdf}
If $0<\vert u\vert\leq q^{b+1/2}$, $b\geq1$ and $x\neq0$, then
    \begin{multline}
        \sum_{n=0}^{\infty}q^{(b+\frac{1}{2})\binom{n}{2}}\rr_{\minus n}(x,y;u|q^{b+\frac{1}{2}})\frac{z^n}{(q;q)_{n}}\\
        =\sum_{i=0}^{\infty}u^{\binom{i}{2}}\frac{(y/x)^{i}}{(\sqrt{q}q^b;\sqrt{q}q^b)_{i}}\sum_{n=0}^{\infty}(q^{-(b+\frac{1}{2})}u)^{\binom{n}{2}}\frac{(-q^{-(b+\frac{1}{2})i}u^iy/x)^{n}}{(\sqrt{q}q^b;\sqrt{q}q^b)_{n}}\\
        \times{}_{0}\phi_{2b}\left(\begin{array}{c}
         -\\
         -\sqrt{q},\mathbf{0}_{2b-1}
    \end{array};\sqrt{q},(-q^{n/2})^{2b+1}z/x\right).
    \end{multline}
\end{theorem}
\begin{proof}
Follow from Theorem \ref{theo_Eboper}. We that
    \begin{align*}
        &\sum_{n=0}^{\infty}q^{(b+\frac{1}{2})\binom{n}{2}}\rr_{\minus n}(x,y;u|q^{b+\frac{1}{2}})\frac{z^n}{(q;q)_{n}}\\
        &\hspace{1cm}=\E(yD_{q^{b+\frac{1}{2}}}|u)\left\{\sum_{n=0}^{\infty}q^{(b+\frac{1}{2})\binom{n}{2}}\frac{(z/x)^n}{(q;q)_{n}}\right\}\\
        &\hspace{1cm}=\E(yD_{q^{b+\frac{1}{2}}}|u)\left\{\mathcal{E}_{q,b}(z/x)\right\}\\
        &\hspace{1cm}=\sum_{i=0}^{\infty}u^{\binom{i}{2}}\frac{(y/x)^{i}}{(\sqrt{q}q^b;\sqrt{q}q^b)_{i}}\sum_{n=0}^{\infty}(q^{-(b+\frac{1}{2})}u)^{\binom{n}{2}}\frac{(-q^{-(b+\frac{1}{2})i}u^iy/x)^{n}}{(\sqrt{q}q^b;\sqrt{q}q^b)_{n}}\\
        &\hspace{4cm}\times{}_{0}\phi_{2b}\left(\begin{array}{c}
         -\\
         -\sqrt{q},\mathbf{0}_{2b-1}
    \end{array};\sqrt{q},(-q^{n/2})^{2b+1}a/x\right).
    \end{align*}
The claim is proved.    
\end{proof}

\subsection{$q^{b}$-series, $b\geq1$}

\begin{lemma}\label{lemma_Rbnder}
For all $n\geq0$, $b\geq2$ and $x\neq0$,
\begin{align}
        D_{q^b}^{n}\{\mathcal{R}_{q,b}(a/x)\}&=(\minus1)^nq^{\minus b\binom{n}{2}}x^{\minus n}\sum_{i=0}^{n}\qbinom{n}{i}_{q^b}(\minus1)^iq^{b\binom{i}{2}}\mathcal{R}_{q,b}(aq^{b(i-n)}/x).
    \end{align}
\end{lemma}
\begin{proof}
Replace $q$ by $q^{b}$ and set $u=q^{b}$ in Lemma \ref{lemma_nder_basic}.    
\end{proof}

\begin{theorem}\label{theo_Rboper}
If $0<\vert u\vert\leq q^b$, $b\geq1$, and $x\neq0$, then
    \begin{multline}
        \E(yD_{q^b}|u)\left\{\mathcal{R}_{q,b}(a/x)\right\}\\
        =\sum_{i=0}^{\infty}u^{\binom{i}{2}}\frac{(y/x)^{i}}{(q^b;q^b)_{i}}\sum_{n=0}^{\infty}(q^{-b}u)^{\binom{n}{2}}\frac{(-q^{-bi}u^iy/x)^{n}}{(q^b;q^b)_{n}}{}_{0}\phi_{b-1}\left(\begin{array}{c}
         -\\
         \mathbf{0}_{b-1}
    \end{array};q,(-1)^{b-1}aq^{\minus bn}/x\right).
    \end{multline}
\end{theorem}
\begin{proof}
Use Lemma \ref{lemma_Rbnder} with $b\geq1$ and then the representation of the function $\mathcal{R}_{q,b}(z)$. 
    \begin{align*}
        &\E(yD_{q^b}|u)\left\{\mathcal{R}_{q,b}(a/x)\right\}\\
        &=\sum_{n=0}^{\infty}u^{\binom{n}{2}}\frac{y^n}{(q^b;q^b)_{n}}D_{q^b}^{n}\{\mathcal{R}_{q,b}(a/x)\}\\
        &=\sum_{n=0}^{\infty}(q^{-b}u)^{\binom{n}{2}}\frac{(\minus y/x)^n}{(q^b;q^b)_{n}}\sum_{i=0}^{n}\qbinom{n}{i}_{q^b}(\minus1)^iq^{b\binom{i}{2}}\mathcal{R}_{q,b}(aq^{b(i-n)}/x)\\
        &=\sum_{i=0}^{\infty}u^{\binom{i}{2}}\frac{(y/x)^{i}}{(q^b;q^b)_{i}}\sum_{n=0}^{\infty}(q^{-b}u)^{\binom{n}{2}}\frac{(-q^{-bi}u^iy/x)^{n}}{(q^b;q^b)_{n}}\mathcal{R}_{q,b}(aq^{\minus bn}/x).
    \end{align*}
The claim is proved.    
\end{proof}

\begin{theorem}\label{theo_Rbgdf}
If $0<\vert u\vert\leq q^b$, $b\geq1$,and $x\neq0$, then
    \begin{multline}
        \sum_{n=0}^{\infty}q^{b\binom{n}{2}}\rr_{\minus n}(x,y;u|q^b)\frac{z^n}{(q;q)_{n}}\\
        =\sum_{i=0}^{\infty}u^{\binom{i}{2}}\frac{(y/x)^{i}}{(q^b;q^b)_{i}}\sum_{n=0}^{\infty}(q^{-b}u)^{\binom{n}{2}}\frac{(-q^{-bi}u^iy/x)^{n}}{(q^b;q^b)_{n}}{}_{0}\phi_{b-1}\left(\begin{array}{c}
         -\\
         \mathbf{0}_{b-1}
    \end{array};q,(-1)^{b-1}aq^{\minus bn}/x\right).
    \end{multline}
\end{theorem}
\begin{proof}
Use Theorem \ref{theo_Rboper}. Hence
    \begin{align*}
        &\sum_{n=0}^{\infty}q^{b\binom{n}{2}}\rr_{\minus n}(x,y;u|q^b)\frac{z^n}{(q;q)_{n}}\\
        &\hspace{1cm}=\E(yD_{q^b}|u)\left\{\sum_{n=0}^{\infty}q^{b\binom{n}{2}}\frac{(z/x)^n}{(q;q)_{n}}\right\}\\
        &\hspace{1cm}=\E(yD_{q^b}|u)\left\{\mathcal{R}_{q,b}(z/x)\right\}\\
        &\hspace{1cm}=\sum_{i=0}^{\infty}u^{\binom{i}{2}}\frac{(y/x)^{i}}{(q^b;q^b)_{i}}\sum_{n=0}^{\infty}(q^{-b}u)^{\binom{n}{2}}\frac{(-q^{-bi}u^iy/x)^{n}}{(q^b;q^b)_{n}}\mathcal{R}_{q,b}(zq^{\minus bn}/x).
    \end{align*}
The claim is reached.    
\end{proof}

\begin{theorem}\label{theo_Roper}
If $0<\vert u\vert\leq q^5$ and $x\neq0$, then
    \begin{multline}
        \E(yD_{q^5}|u)\left\{\mathcal{R}_{q}(a/x)\right\}
        =\frac{(aq^5/x,(aq/x)^2,(aq/x)^2q;q^5)_{\infty}}{(aq/x;q)_{\infty}}\sum_{i=0}^{\infty}\frac{u^{\binom{i}{2}}(y/x)^{i}}{(q^5;q^5)_{i}}\\
        \hspace{1cm}\times\sum_{n=0}^{\infty}\left(\frac{u}{q^5}\right)^{\binom{n}{2}}\frac{(x/a;q)_{n}((xq/a)^{2}q,(xq/a)^2;q^5)_{2n}}{(x/a;q)_{5n}(q^5;q^5)_{n}}(-q^{-5(i-1)}u^ia^4y/x^5)^{n}\\
        \hspace{2cm}\times
        {}_{3}\phi_{2}\left(\begin{array}{c}
         aq^{-5n-1}/x,aq^{-5n}/x,aq^{-5n+1}/x\\
         a^2q^{-10n+2}/x^2,a^2q^{-10n+3}/x^2
    \end{array};q^5,aq^{-5(n-1)}/x\right).
    \end{multline}
\end{theorem}
\begin{proof}
Follows from Lemma \ref{lemma_Rbnder} with $b=5$,
    \begin{align*}
        &\E(yD_{q^5}|u)\left\{\mathcal{R}_{q}(a/x)\right\}\\
        &=\sum_{n=0}^{\infty}u^{\binom{n}{2}}\frac{y^n}{(q^5;q^5)_{n}}D_{q^5}^{n}\{\mathcal{R}_{q}(a/x)\}\\
        &=\sum_{n=0}^{\infty}(q^{-5}u)^{\binom{n}{2}}\frac{(\minus y/x)^n}{(q^5;q^5)_{n}}\sum_{i=0}^{n}\qbinom{n}{i}_{q^5}(\minus1)^iq^{5\binom{i}{2}}\mathcal{R}_{q}(aq^{5(i-n)}/x)\\
        &=\sum_{i=0}^{\infty}u^{\binom{i}{2}}\frac{(y/x)^{i}}{(q^5;q^5)_{i}}\sum_{n=0}^{\infty}(q^{-5}u)^{\binom{n}{2}}\frac{(-q^{-5i}u^ix^{\minus1}y)^{n}}{(q^5;q^5)_{n}}\mathcal{R}_{q}(aq^{\minus5n}/x).
    \end{align*}
Now, by applying the Lemma \ref{lemma_prop_exp4},    
    \begin{align*}
    &\E(yD_{q^5}|u)\left\{\mathcal{R}_{q}(a/x)\right\}\\
        &=\frac{(aq^5/x,(aq/x)^2,(aq/x)^2q;q^5)_{\infty}}{(aq/x;q)_{\infty}}\sum_{i=0}^{\infty}\frac{u^{\binom{i}{2}}(y/x)^{i}}{(q^5;q^5)_{i}}\\
        &\hspace{1cm}\times\sum_{n=0}^{\infty}\left(\frac{u}{q^5}\right)^{\binom{n}{2}}\frac{(x/a;q)_{n}((xq/a)^{2}q,(xq/a)^2;q^5)_{2n}}{(x/a;q)_{5n}(q^5;q^5)_{n}}(-q^{-5(i-1)}u^ia^4y/x^5)^{n}\\
        &\hspace{2cm}\times
        {}_{3}\phi_{2}\left(\begin{array}{c}
         aq^{-5n-1}/x,aq^{-5n}/x,aq^{-5n+1}/x\\
         a^2q^{-10n+2}/x^2,a^2q^{-10n+3}/x^2
    \end{array};q^5,aq^{-5(n-1)}/x\right).
    \end{align*}
The claim is proved.    
\end{proof}

\begin{theorem}\label{theo_Rgdf}
If $0<\vert u\vert\leq q^5$ and $x\neq0$, then
    \begin{multline}
        \sum_{n=0}^{\infty}q^{n^2}\rr_{\minus n}(x,y;u|q^5)\frac{z^n}{(q;q)_{n}}
        =\frac{(zq^5/x,(zq/x)^2,(zq/x)^2q;q^5)_{\infty}}{(zq/x;q)_{\infty}}\sum_{i=0}^{\infty}\frac{u^{\binom{i}{2}}(y/x)^{i}}{(q^5;q^5)_{i}}\\
        \hspace{1cm}\times\sum_{n=0}^{\infty}\left(\frac{u}{q^5}\right)^{\binom{n}{2}}\frac{(x/z;q)_{n}((xq/z)^{2}q,(xq/z)^2;q^5)_{2n}}{(x/z;q)_{5n}(q^5;q^5)_{n}}(-q^{-5(i-1)}u^iz^4y/x^5)^{n}\\
        \hspace{2cm}\times
        {}_{3}\phi_{2}\left(\begin{array}{c}
         zq^{-5n-1}/x,zq^{-5n}/x,zq^{-5n+1}/x\\
         z^2q^{-10n+2}/x^2,z^2q^{-10n+3}/x^2
    \end{array};q^5,zq^{-5(n-1)}/x\right).
    \end{multline}
\end{theorem}
\begin{proof}
Use Theorem \ref{theo_Roper} for noting that
    \begin{align*}
        &\sum_{n=0}^{\infty}q^{n^2}\rr_{\minus n}(x,y;u|q^5)\frac{z^n}{(q;q)_{n}}\\
        &\hspace{1cm}=\E(yD_{q^5}|u)\left\{\sum_{n=0}^{\infty}q^{n^2}\frac{(z/x)^n}{(q;q)_{n}}\right\}\\
        &\hspace{1cm}=\E(yD_{q^5}|u)\left\{\mathcal{R}_{q}(z/x)\right\}\\
        &\hspace{1cm}=\frac{(zq^5/x,(zq/x)^2,(zq/x)^2q;q^5)_{\infty}}{(zq/x;q)_{\infty}}\sum_{i=0}^{\infty}\frac{u^{\binom{i}{2}}(y/x)^{i}}{(q^5;q^5)_{i}}\\
        &\hspace{2cm}\times\sum_{n=0}^{\infty}\left(\frac{u}{q^5}\right)^{\binom{n}{2}}\frac{(x/z;q)_{n}((xq/z)^{2}q,(xq/z)^2;q^5)_{2n}}{(x/z;q)_{5n}(q^5;q^5)_{n}}(-q^{-5(i-1)}u^iz^4y/x^5)^{n}\\
        &\hspace{2cm}\times
        {}_{3}\phi_{2}\left(\begin{array}{c}
         zq^{-5n-1}/x,zq^{-5n}/x,zq^{-5n+1}/x\\
         z^2q^{-10n+2}/x^2,z^2q^{-10n+3}/x^2
    \end{array};q^5,zq^{-5(n-1)}/x\right).
    \end{align*}
\end{proof}

\section{Some special cases}

\subsection{$q$-Series involving the function of Cauchy $\P_{\minus n}(x,y;q)$}

\begin{corollary}
If $u=q$ in Theorem \ref{theo_operE1}, then
    \begin{equation}
        \E(yD_{q}|q)\left\{\frac{1}{(a/x;q)_{\infty}}\right\}=\frac{1}{(a/x;q)_{\infty}}\sum_{i=0}^{\infty}q^{\binom{i}{2}}\frac{(y/x)^{i}}{(q;q)_{i}}{}_{1}\phi_{1}\left(
        \begin{array}{c}
             0\\
             qx/a
        \end{array};
        q,q^{i+1}y/a
        \right)
    \end{equation}
for all $0<\vert a/x\vert<1$.    
\end{corollary}

\begin{corollary}\label{coro_ncgf}
If $u=q$ and $y\mapsto-y$ in Theorem \ref{theo_gdfE1}, then for $0<\vert z/x\vert<1$, 
    \begin{align}
        \sum_{n=0}^{\infty}\P_{\minus n}(x,y;q)\frac{z^n}{(q;q)_{n}}
        =\frac{1}{(z/x;q)_{\infty}}\sum_{i=0}^{\infty}q^{\binom{i}{2}}\frac{(-y/x)^{i}}{(q;q)_{i}}{}_{1}\phi_{1}\left(
        \begin{array}{c}
             0\\
             qx/z
        \end{array};
        q,-q^{i+1}y/z
        \right).
    \end{align}
\end{corollary}
\begin{corollary}
Suppose that $0<\vert z\vert<1$ and $x\neq0$. Then the transformation formula holds true
    \begin{equation}
        {}_{1}\phi_{1}\left(\begin{array}{c}
         0  \\
         q/x
    \end{array};q,qz/x\right)=\frac{1}{(z;q)_{\infty}}\sum_{i=0}^{\infty}q^{\binom{i}{2}}\frac{(-x)^{i}}{(q;q)_{i}}{}_{1}\phi_{1}\left(\begin{array}{c}
         0  \\
         q/z
    \end{array};q,\minus q^{i+1}x/z\right).    
    \end{equation}
\end{corollary}
\begin{proof}
If $x=1$ and $y=x$ in Corollary \ref{coro_ncgf} and as
\begin{equation}\label{iden10}
    \P_{\minus n}(1,x;q)=(x;q)_{\minus n}=q^{\binom{n}{2}}\frac{(-q/x)^n}{(q/x;q)_{n}},
\end{equation}
then
\begin{equation}
    \sum_{n=0}^{\infty}q^{\binom{n}{2}}\frac{(\minus qz/x)^n}{(q/x;q)_{n}(q;q)_{n}}=\frac{1}{(z;q)_{\infty}}\sum_{i=0}^{\infty}q^{\binom{i}{2}}\frac{(-x)^{i}}{(q;q)_{i}}{}_{1}\phi_{1}\left(\begin{array}{c}
         0  \\
         q/z
    \end{array};q,\minus q^{i+1}x/z\right).
\end{equation}
Hence, we obtain the transformation formula
\begin{equation}
    {}_{1}\phi_{1}\left(\begin{array}{c}
         0  \\
         q/x
    \end{array};q,qz/x\right)=\frac{1}{(z;q)_{\infty}}\sum_{i=0}^{\infty}q^{\binom{i}{2}}\frac{(-x)^{i}}{(q;q)_{i}}{}_{1}\phi_{1}\left(\begin{array}{c}
         0  \\
         q/z
    \end{array};q,\minus q^{i+1}x/z\right).
\end{equation}
\end{proof}

\begin{corollary}\label{coro_Eoper}
If $u=\sqrt{q}$ in Theorem \ref{theo_Eoper} and $x\neq0$, then
    \begin{multline}
        \E(yD_{\sqrt{q}}|\sqrt{q})\left\{\mathcal{E}_{q}(z/x)\right\}\\
        =(-y/x;\sqrt{q})_{\infty}\sum_{n=0}^{\infty}\frac{(-y/x)^{n}}{(\sqrt{q};\sqrt{q})_{n}}{}_{1}\phi_{1}\left(\begin{array}{c}
         0\\
         -\sqrt{q}
    \end{array};\sqrt{q},\minus q^{-n}z/x\right).
    \end{multline}
\end{corollary}

\begin{corollary}\label{coro_Egdf}
If $u=\sqrt{q}$ in Theorem \ref{theo_Egdf} and $x\neq0$, then
    \begin{multline}
        \sum_{n=0}^{\infty}q^{\frac{1}{2}\binom{n}{2}}\P_{\minus n}(x,y;\sqrt{q})\frac{z^n}{(q;q)_{n}}\\
        =(-y/x;\sqrt{q})_{\infty}\sum_{n=0}^{\infty}\frac{(-y/x)^{n}}{(\sqrt{q};\sqrt{q})_{n}}{}_{1}\phi_{1}\left(\begin{array}{c}
         0\\
         -\sqrt{q}
    \end{array};\sqrt{q},\minus q^{-n}z/x\right).
    \end{multline}
\end{corollary}
If $x=1$ and $y=x$ in Corollary \ref{coro_Egdf}, then
\begin{corollary}
Suppose that $x\neq0$. Then we have the transformation formula
    \begin{equation}
        {}_{1}\phi_{2}\left(\begin{array}{c}
         0\\
         q/x,-q
    \end{array};q,-qz/x\right)
    =(-x;q)_{\infty}\sum_{n=0}^{\infty}\frac{(-x)^{n}}{(q;q)_{n}}{}_{1}\phi_{1}\left(\begin{array}{c}
         0\\
         -q
    \end{array};q,-q^{-2n}z\right).
    \end{equation}
\end{corollary}
\begin{proof}
On the one side,
    \begin{align*}
        &\sum_{n=0}^{\infty}q^{\frac{1}{2}\binom{n}{2}}(x;\sqrt{q})_{-n}\frac{z^n}{(-\sqrt{q};\sqrt{q})_{n}(\sqrt{q};\sqrt{q})_{n}}\\
        &\hspace{1cm}=\sum_{n=0}^{\infty}q^{\binom{n}{2}}\frac{(-\sqrt{q}z/x)^n}{(\sqrt{q}/x;\sqrt{q})_{n}(-\sqrt{q};\sqrt{q})_{n}(\sqrt{q};\sqrt{q})_{n}}\\
        &\hspace{1cm}={}_{1}\phi_{2}\left(\begin{array}{c}
         0\\
         \sqrt{q}/x,-\sqrt{q}
    \end{array};\sqrt{q},-\sqrt{q}z/x\right).
    \end{align*}
On the other hand,
\begin{align*}
    &\sum_{n=0}^{\infty}q^{\frac{1}{2}\binom{n}{2}}(x;\sqrt{q})_{-n}\frac{z^n}{(-\sqrt{q};\sqrt{q})_{n}(\sqrt{q};\sqrt{q})_{n}}\\
    &=(-x;\sqrt{q})_{\infty}\sum_{n=0}^{\infty}\frac{(-x)^{n}}{(\sqrt{q};\sqrt{q})_{n}}{}_{1}\phi_{1}\left(\begin{array}{c}
         0\\
         -\sqrt{q}
    \end{array};\sqrt{q},\minus q^{-n}z\right).
\end{align*}
Finally, replace $\sqrt{q}$ by $q$.
\end{proof}

\begin{corollary}
If $u=q^{b+1/2}$ in Theorem \ref{theo_Eboper}, with $b\geq1$, and $x\neq0$ then
    \begin{multline}
        \E(yD_{q^{b+1/2}}|q^{b+1/2})\left\{\mathcal{E}_{q,b}(a/x)\right\}\\
        =(-y/x;q^{b+1/2})_{\infty}\sum_{n=0}^{\infty}\frac{(-y/x)^{n}}{(\sqrt{q}q^b;\sqrt{q}q^b)_{n}}{}_{0}\phi_{2b}\left(\begin{array}{c}
         -\\
         -\sqrt{q},\mathbf{0}_{2b-1}
    \end{array};\sqrt{q},(-q^{n/2})^{2b+1}a/x\right).
    \end{multline}
\end{corollary}

\begin{corollary}\label{coro_Ebgdf}
If $u=q^{b+1/2}$ in Theorem \ref{theo_Ebgdf}, $b\geq1$, and $x\neq0$, then
    \begin{multline}
        \sum_{n=0}^{\infty}q^{(b+\frac{1}{2})\binom{n}{2}}\P_{\minus n}(x,y;q^{b+\frac{1}{2}})\frac{z^n}{(q;q)_{n}}\\
        =(-y/x;q^{b+\frac{1}{2}})_{\infty}\sum_{n=0}^{\infty}\frac{(-y/x)^{n}}{(\sqrt{q}q^b;\sqrt{q}q^b)_{n}}
        {}_{0}\phi_{2b}\left(\begin{array}{c}
         -\\
         -\sqrt{q},\mathbf{0}_{2b-1}
    \end{array};\sqrt{q},(-q^{n/2})^{2b+1}z/x\right).
    \end{multline}
\end{corollary}
If $x=1$ and $y=x$ in Corollary \ref{coro_Ebgdf}, then
\begin{corollary}
For all $b\geq1$, and $x\neq0$,
    \begin{multline}
        \sum_{n=0}^{\infty}\frac{q^{2(2b+1)\binom{n}{2}}(-q^{2b+1}z/x)^n}{(q^{2b+1}/x;q^{2b+1})_{n}(-q;q)_{n}(q;q)_{n}}\\
        =(-x;q^{2b+1})_{\infty}\sum_{n=0}^{\infty}\frac{(-x)^{n}}{(q^{2b+1};q^{2b+1})_{n}}
        {}_{0}\phi_{2b}\left(\begin{array}{c}
         -\\
         -q,\mathbf{0}_{2b-1}
    \end{array};q,(-q^{n})^{2b+1}z\right).
    \end{multline}
\end{corollary}
\begin{proof}
On the one side,
\begin{align*}
    &\sum_{n=0}^{\infty}q^{(b+\frac{1}{2})\binom{n}{2}}(x;q^{b+\frac{1}{2}})_{-n}\frac{z^n}{(-\sqrt{q};\sqrt{q})_{n}(\sqrt{q};\sqrt{q})_{n}}\\
    &\hspace{1cm}=\sum_{n=0}^{\infty}q^{(2b+1)\binom{n}{2}}\frac{(-q^{(b+\frac{1}{2})}z/x)^n}{(q^{(b+\frac{1}{2})}/x;q^{b}\sqrt{q})_{n}(-\sqrt{q};\sqrt{q})_{n}(\sqrt{q};\sqrt{q})_{n}}.
\end{align*}
On the other hand,
\begin{align*}
    &\sum_{n=0}^{\infty}q^{(b+\frac{1}{2})\binom{n}{2}}(x;q^{b+\frac{1}{2}})_{-n}\frac{z^n}{(-\sqrt{q};\sqrt{q})_{n}(\sqrt{q};\sqrt{q})_{n}}\\
    &\hspace{1cm}=(-x;q^{b+\frac{1}{2}})_{\infty}\sum_{n=0}^{\infty}\frac{(-x)^{n}}{(\sqrt{q}q^b;\sqrt{q}q^b)_{n}}
        {}_{0}\phi_{2b}\left(\begin{array}{c}
         -\\
         -\sqrt{q},\mathbf{0}_{2b-1}
    \end{array};\sqrt{q},(-q^{n/2})^{2b+1}z\right).
\end{align*}
Finally, replace $\sqrt{q}$ with $q$.
\end{proof}

\begin{corollary}\label{coro_Rboper}
If $u=q^b$ in Theorem \ref{theo_Rboper} and $x\neq0$, then
    \begin{multline}
        \E(yD_{q^b}|q^{b})\left\{\mathcal{R}_{q,b}(z/x)\right\}\\
        =(y/x;q^b)_{\infty}\sum_{n=0}^{\infty}\frac{(-y/x)^{n}}{(q^b;q^b)_{n}}{}_{0}\phi_{b-1}\left(\begin{array}{c}
         -\\
         \mathbf{0}_{b-1}
    \end{array};q,(-1)^{b-1}zq^{\minus bn}/x\right).
    \end{multline}
\end{corollary}

\begin{corollary}\label{coro_Rbgdf}
If $u=q^b$ in Theorem \ref{theo_Rbgdf}, then
    \begin{multline}
        \sum_{n=0}^{\infty}q^{b\binom{n}{2}}\P_{\minus n}(x,y;q^b)\frac{z^n}{(q;q)_{n}}\\
        =(y/x;q^b)\sum_{n=0}^{\infty}\frac{(-y/x)^{n}}{(q^b;q^b)_{n}}{}_{0}\phi_{b-1}\left(\begin{array}{c}
         -\\
         \mathbf{0}_{b-1}
    \end{array};q,(-1)^{b-1}aq^{\minus bn}/x\right).
    \end{multline}
\end{corollary}

If $x=1$ and $y=x$ in Corollary \ref{coro_Rbgdf}, then
\begin{corollary}
For all $b\geq1$ and $x\neq0$,
    \begin{align}
        \sum_{n=0}^{\infty}q^{2b\binom{n}{2}}\frac{(-q^bz/x)^n}{(q^b/x;q^b)_{n}(q;q)_{n}}
        =(x;q^b)_{\infty}\sum_{n=0}^{\infty}\frac{(-x)^{n}}{(q^b;q^b)_{n}}{}_{0}\phi_{b-1}\left(\begin{array}{c}
         -\\
         \mathbf{0}_{b-1}
    \end{array};q,(-1)^{b-1}q^{\minus bn}z\right).
    \end{align}
\end{corollary}
\begin{proof}
On the one side,
    \begin{align*}
        \sum_{n=0}^{\infty}q^{b\binom{n}{2}}(x;q^b)_{\minus n}\frac{z^n}{(q;q)_{n}}&=\sum_{n=0}^{\infty}q^{2b\binom{n}{2}}\frac{(-q^bz/x)^n}{(q^b/x;q^b)_{n}(q;q)_{n}}.
    \end{align*}
On the other hand,
\begin{align*}\sum_{n=0}^{\infty}q^{b\binom{n}{2}}(x;q^b)_{\minus n}\frac{z^n}{(q;q)_{n}}&=
    (x;q^b)_{\infty}\sum_{n=0}^{\infty}\frac{(-x)^{n}}{(q^b;q^b)_{n}}\mathcal{R}_{q,b}(zq^{\minus bn}).
\end{align*}
\end{proof}

\begin{corollary}
If $u=q^5$ and $x\neq0$, then
    \begin{multline}
        \sum_{n=0}^{\infty}q^{n^2}\P_{\minus n}(x,y;q^5)\frac{z^n}{(q;q)_{n}}=\frac{(y/x,zq^5/x,(zq/x)^2,(zq/x)^2q;q^5)_{\infty}}{(zq/x;q)_{\infty}}\\
        \hspace{1cm}\times\sum_{n=0}^{\infty}\frac{(x/z;q)_{n}((xq/z)^{2}q,(xq/z)^2;q^5)_{2n}}{(x/z;q)_{5n}(q^5;q^5)_{n}}(q^{5}z^4y/x^5)^{n}\\
        \hspace{2cm}\times
        {}_{3}\phi_{2}\left(\begin{array}{c}
         zq^{-5n-1}/x,zq^{-5n}/x,zq^{-5n+1}/x\\
         z^2q^{-10n+2}/x^2,z^2q^{-10n+3}/x^2
    \end{array};q^5,zq^{-5(n-1)}/x\right),
    \end{multline}
for all $0<\vert z/x\vert<1$.    
\end{corollary}
If $x=1$ and $y=x$, then
\begin{align}
        \sum_{n=0}^{\infty}q^{\frac{7n^2-5n}{2}}\frac{(-q^5z/x)^n}{(q^5/x;q^5)_{n}(q;q)_{n}}&=\frac{(x,zq^5,(zq)^2,(zq)^2q;q^5)_{\infty}}{(zq;q)_{\infty}}\nonumber\\
        &\hspace{0.5cm}\times\sum_{n=0}^{\infty}\frac{(1/z;q)_{n}((q/z)^{2}q,(q/z)^2;q^5)_{2n}}{(1/z;q)_{5n}(q^5;q^5)_{n}}(q^{5}z^4x)^{n}\nonumber\\
        &\hspace{0.5cm}\times
        {}_{3}\phi_{2}\left(\begin{array}{c}
         zq^{-5n-1},zq^{-5n},zq^{-5n+1}\\
         z^2q^{-10n+2},z^2q^{-10n+3}
    \end{array};q^5,zq^{-5(n-1)}\right).
    \end{align}

\subsection{$q$-Series involving the function of Stieltjes-Wigert $\s_{\minus n}(x,y;q)$}

\begin{corollary}
If $u=q^2$ in Theorem \ref{theo_operE1} and $x\neq0$, $a\neq0$, then
    \begin{equation}
        \E(qyD_{q}|q^2)\left\{\frac{1}{(a/x;q)_{\infty}}\right\}=\frac{1}{(a/x;q)_{\infty}}\sum_{i=0}^{\infty}q^{i^2}\frac{(y/x)^{i}}{(q;q)_{i}}{}_{0}\phi_{1}\left(
        \begin{array}{c}
             - \\
             qx/a
        \end{array};
        q,q^{2i+2}y/a
        \right),
    \end{equation}
for all $\vert a/x\vert<1$.
\end{corollary}

\begin{corollary}
If $u=q^2$ in Theorem \ref{theo_gdfE1} then
    \begin{align}
        \sum_{n=0}^{\infty}\s_{\minus n}(x,y;q)\frac{z^n}{(q;q)_{n}}=\frac{1}{(z/x;q)_{\infty}}\sum_{i=0}^{\infty}q^{i^2}\frac{(y/x)^{i}}{(q;q)_{i}}{}_{0}\phi_{1}\left(
        \begin{array}{c}
             - \\
             qx/z
        \end{array};
        q,q^{2i+2}y/z
        \right)
    \end{align}
for all $0<\vert z/x\vert<1$.    
\end{corollary}

\begin{corollary}
If $u=q^2$ in Theorem \ref{theo_operE2}, then
    \begin{multline}
        \E(qyD_{q}|q^2)\left\{(z/x;q)_{\infty}\right\}\\
        =\frac{(z/x,qy/x;q)_{\infty}}{(zy/x^2;q)_{\infty}}{}_{4}\phi_{5}\left(
        \begin{array}{c}
             \sqrt{zy/x^2},-\sqrt{zy/x^2},\sqrt{zqy/x^2},-\sqrt{zqy/x^2} \\
             \sqrt{qy/x},\sqrt{qy/x},\sqrt{q^2y/x},\sqrt{q^2y/x},0
        \end{array};
        q,qy/x
        \right)
    \end{multline}
for all $0<\vert zy/x^2\vert<1$.    
\end{corollary}

\begin{corollary}
If $u=q^2$ in Theorem \ref{theo_gdfE2}, then
    \begin{multline}
        \sum_{n=0}^{\infty}(-1)^nq^{\binom{n}{2}}\s_{\minus n}(x,y;q)\frac{z^n}{(q;q)_{n}}\\
        =\frac{(z/x,qy/x;q)_{\infty}}{(zy/x^2;q)_{\infty}}{}_{4}\phi_{5}\left(
        \begin{array}{c}
             \sqrt{zy/x^2},-\sqrt{zy/x^2},\sqrt{zqy/x^2},-\sqrt{zqy/x^2} \\
             \sqrt{qy/x},\sqrt{qy/x},\sqrt{q^2y/x},\sqrt{q^2y/x},0
        \end{array};
        q,qy/x
        \right)
    \end{multline}
for all $0<\vert zy/x^2\vert<1$.    
\end{corollary}

\subsection{$q$-Series involving the function $\s_{\minus n,b}(x,y;q)$, $b\geq3$}

\begin{corollary}
If $u=q^b$ in Theorem \ref{theo_operE1}, then
    \begin{multline}
        \E(yD_{q}|q^b)\left\{\frac{1}{(a/x;q)_{\infty}}\right\}\\
        =\frac{1}{(a/x;q)_{\infty}}\sum_{i=0}^{\infty}q^{b\binom{i}{2}}\frac{(y/x)^{i}}{(q;q)_{i}}{}_{2}\phi_{1+b}\left(
        \begin{array}{c}
             0,0 \\
             qx/a,\mathbf{0}_{b}
        \end{array};
        q,(-1)^bq^{bi+1}y/a
        \right),
    \end{multline}
for all $0<\vert a/x\vert<1$.    
\end{corollary}

\begin{corollary}
If $u=q^{b}$ in Theorem \ref{theo_gdfE1}, then
    \begin{multline}
        \sum_{n=0}^{\infty}\s_{\minus n,b}(x,y;q)\frac{z^n}{(q;q)_{n}}\\
        =\frac{1}{(z/x;q)_{\infty}}\sum_{i=0}^{\infty}q^{b\binom{i}{2}}\frac{(y/x)^{i}}{(q;q)_{i}}{}_{2}\phi_{1+b}\left(
        \begin{array}{c}
             0,0 \\
             qx/z,\mathbf{0}_{b}
        \end{array};
        q,(-1)^bq^{bi+1}y/z
        \right),
    \end{multline}
for all $0<\vert z/x\vert<1$.    
\end{corollary}

\begin{corollary}
If $b\geq3$ in Theorem \ref{theo_operE2} and $x\neq0$, $z\neq0$, then
    \begin{multline}
        \E(yD_{q}|q^b)\left\{(z/x;q)_{\infty}\right\}\\
        =(z/x;q)_{\infty}\sum_{i=0}^{\infty}q^{b\binom{i}{2}}\frac{(y/x)^{i}}{(q;q)_{i}}{}_{1}\phi_{b-2}\left(
        \begin{array}{c}
             qx/z \\
             \mathbf{0}_{b-2}
        \end{array};
        q,(-1)^{b-2}q^{bi-1}zy/x^2
        \right).
    \end{multline}
\end{corollary}

\begin{corollary}
If $b\geq3$ in Theorem \ref{theo_gdfE2} and $x\neq0$, $z\neq0$, then
    \begin{multline}
        \sum_{n=0}^{\infty}(-1)^nq^{\binom{n}{2}}\s_{\minus n,b}(x,y;q)\frac{z^n}{(q;q)_{n}}\\
        =(z/x;q)_{\infty}\sum_{i=0}^{\infty}q^{b\binom{i}{2}}\frac{(y/x)^{i}}{(q;q)_{i}}{}_{1}\phi_{b-2}\left(
        \begin{array}{c}
             qx/z \\
             \mathbf{0}_{b-2}
        \end{array};
        q,(-1)^{b-2}q^{bi-1}zy/x^2
        \right).
    \end{multline}
\end{corollary}

\subsection{$q$-Series involving the function $\E_{\minus n,b}(x,y;q)$, $b\geq1$}

\begin{corollary}
If $u=q^{b+1/2}$ in Theorem \ref{theo_operE1}, then
    \begin{multline}
        \E(yD_{q}|\sqrt{q})\left\{\frac{1}{(z/x;q)_{\infty}}\right\}
        =\frac{1}{(z/x;q)_{\infty}}\sum_{i=0}^{\infty}q^{(b+\frac{1}{2})\binom{i}{2}}\frac{(y/x)^{i}}{(q;q)_{i}}\\\times{}_{4}\phi_{2b+4}\left(
        \begin{array}{c}
             0,0,0,0 \\
             \sqrt{qx/z},-\sqrt{qx/z},-\sqrt{q},\mathbf{0}_{2b+1}
        \end{array};
        \sqrt{q},-q^{i(b+\frac{1}{2})+1}y/z
        \right),
    \end{multline}
for all $0<\vert z/x\vert<1$.    
\end{corollary}

\begin{corollary}
If $u=q^{b+1/2}$ in Theorem \ref{theo_gdfE1}, then
    \begin{multline}
        \sum_{n=0}^{\infty}\E_{\minus n,b}(x,y;q)\frac{z^n}{(q;q)_{n}}
        =\frac{1}{(z/x;q)_{\infty}}\sum_{i=0}^{\infty}q^{(b+\frac{1}{2})\binom{i}{2}}\frac{(y/x)^{i}}{(q;q)_{i}}\\\times{}_{4}\phi_{2b+4}\left(
        \begin{array}{c}
             0,0,0,0 \\
             \sqrt{qx/z},-\sqrt{qx/z},-\sqrt{q},\mathbf{0}_{2b+1}
        \end{array};
        \sqrt{q},-q^{i(b+\frac{1}{2})+1}y/z
        \right),
    \end{multline}
for all $0<\vert z/x\vert<1$.    
\end{corollary}


\begin{thebibliography}{99}

\bibitem{orozco}
R. Orozco, 
Deformed homogeneous polynomials and the deformed $q$-exponential operator, 
Commun. Math. \textbf{34}, 2, 7
(2024).

\bibitem{chen1}
W.Y.C. Chen, Z.G. Liu,
Parameter augmenting for basic hypergeometric series, I,
Mathematical Essays in Honour of Gian-Carlo Rota, Eds., B. E. Sagan and R.P. Stanley, Birkh\"auser, Boston, 1998, pp. 111--129.

\bibitem{chen2}
W.Y.C. Chen, Z.G. Liu,
Parameter augmenting for basic hypergeometric series, II,
J. Combin. Theory, Ser. A
\textbf{80} (1997) 175--195.

\bibitem{gasper}
G. Gasper, M. Rahman,
Basic Hypergeometric Series, $2^{nd}$ ed., 
Cambridge University Press, Cambridge, MA, 1990.

\bibitem{saad}
H. L. Saad and M.A. Abdul, 
The $q$-Exponential Operator and Generalized Rogers-Szego Polynomials, 
Journal of Advances in Mathematics, 8(1), (2014) 1440–-1455. 
https://doi.org/10.24297/jam.v8i1.6912

\bibitem{saad2}
H. L. Saad, and A. A. Sukhi, 
The $q$-exponential operator, 
Applied Mathematical Sciences, \textbf{7}, 6369--6380, (2013)
http://dx.doi.org/10.12988/ams.2013.24287.




\end{thebibliography}
\end{document}